\documentclass{article}

\advance\textheight -\headheight
\advance\textheight -\headsep
\oddsidemargin\paperwidth
\advance\oddsidemargin -\textwidth
\divide\oddsidemargin2 
\ifdim\oddsidemargin<.5truein \oddsidemargin.5truein \fi
\advance\oddsidemargin -1truein
\evensidemargin\oddsidemargin
\topmargin\paperheight \advance\topmargin -\textheight
\advance\topmargin -\headheight \advance\topmargin -\headsep
\divide\topmargin2
\ifdim\topmargin<.5truein \topmargin.5truein \fi
\advance\topmargin -1.4truein\relax

\usepackage{graphicx,amssymb,latexsym,amsfonts,txfonts}
\usepackage{pdfsync,color,tabularx,rotating}
\usepackage[all,cmtip]{xy}
\usepackage{url}
\usepackage{tikz}
\usepackage{amssymb}
\usepackage{comment}
\usepackage{epsfig}
\usepackage{eufrak}

\newtheorem{theo}{Theorem}[section]
\newtheorem{prop}[theo]{Proposition}
\newtheorem{coro}[theo]{Corollary}
\newtheorem{lemm}[theo]{Lemma}

\begin{document}

\title{Permutation groups of prime power degree and $p$-complements}

\author{Gareth A. Jones and Sezgin Sezer}

\maketitle

\begin{abstract}

Extending earlier work of Guralnick and of Cai and Zhang, we classify the almost simple groups which have transitive permutation representations of prime power degree $p^k$, and those which have $p$-complements (stabilisers of order coprime to $p$ in such representations).
We deduce that every primitive permutation group of prime power degree has a regular subgroup, and that any two faithful primitive representations of a group, of the same prime power degree, are equivalent under automorphisms. In general, $p$-complements in a finite group can be inequivalent under automorphisms, or even non-isomorphic. We extend examples of such phenomena due to Buturlakin, Revin and Nesterov by showing that the number of inequivalent classes of complements can be arbitrarily large.
Questions concerning the existence of prime power representations and $p$-complements in groups with socle ${\rm PSL}_d(q)$ are related to some difficult open problems in Number Theory. 

\end{abstract}


\medskip



\section{Introduction}
\label{sec:Intro}

In 1983 Guralnick~\cite{Gur} used the then recently-announced classification of non-abelian finite simple groups to determine their transitive permutation representations of prime power degree. In 2015 Cai and Zhang~\cite{CZ} extended this result to classify the primitive permutation groups of prime power degree. Here, motivated by problems concerning $p$-complements discussed later, we give a further partial extension of these results to transitive representations of prime power degree for almost simple groups, those groups $G$ such that $S\le G\le {\rm Aut}\,S$ for some non-abelian finite simple group $S$. Our first main result is as follows:

\begin{theo}\label{th:mainppower}
The almost simple groups which are transitive permutation groups of prime power degree $n=p^k$ are as follows:
\begin{itemize}
\item[\rm(A)]
\begin{enumerate}
\item  ${\rm A}_n$ and ${\rm S}_n$ for $n=p^k\ge 5$;
\item ${\rm S}_{n/2}$ for $n=2^k\ge 16$;
\end{enumerate}
\item[\rm(B)] 
\begin{enumerate}
\item ${\rm PSL}_2(q)$ and ${\rm PGL}_2(q)$ for any Mersenne prime $q=2^k-1>3$, with $p=2$;
\item ${\rm PSL}_2(q)\rtimes{\rm C}_{2^i}\le{\rm P\Gamma L}_2(q)$ for $q=2^e>2$ and $i=0,\ldots,f$, with $e=2^f$ and $p^k=p=2^e+1$ a Fermat prime;
\item ${\rm PSL}_2(8)$ and ${\rm P\Gamma L}_2(8)$, where $p^k=9$;
\item for $d>2$, any group $G$ such that ${\rm PSL}_d(q)\le G\le {\rm P\Gamma L}_d(q)$ where $n=(q^d-1)/(q-1)$ is a prime $p$;
\item for $d>2$, any group $G$ such that ${\rm PSL}_d(q)\le G\le {\rm P\Gamma L}_d(q)$ where $n=(q^d-1)/(q-1)$ is a proper prime power $p^k$, $k>1$;
\end{enumerate}
\item[\rm(C)] ${\rm PSL}_ 2(11)$ and ${\rm M}_{11}$ for $n=11$, and ${\rm M}_{23}$ for $n=23$;
\item[\rm(D)] $W(E_6)$ and $W(E_6)'\cong{\rm PSp}_4(3)$ for $n=27$.
\end{itemize}
There are two such permutation representations for $G$ in {\rm(B4)} and {\rm(B5)}, and for ${\rm PSL}_2(11)$ in {\rm(C)}, transposed by an outer automorphism of $G$ in each case. In all other cases each group has just one such representation. All of these groups except those in {\rm(A2)} are primitive, and of these, all are doubly transitive except those in {\rm(D)}, which have rank $3$.
\end{theo}

In (A1), the alternating and symmetric groups have their natural representations (since $6$ is not a prime power, their additional representations of that degree do not arise); these groups are solvable for $n<5$, so they are excluded. The action in (A2) is the product of the natural action of degree $n/2$ and the signature action on $\{\pm 1\}$. In (B) $q$ is a prime power, and the various groups act naturally on the points (or the points and the hyperplanes if $d>2$) of the associated projective geometry ${\mathbb P}^{d-1}({\mathbb F}_q)$, which has ${\rm P\Gamma L}_d(q)$ as its automorphism group. In (B1), for $q$ to be prime it is necessary that $k$ should be prime; at the time of writing, $51$ Mersenne primes are known~\cite{GIMPS}, and it is conjectured that there exist infinitely many; here $q=3$ is excluded since the corresponding groups are solvable. In (B2) five Fermat primes are currently known, for $f=0, 1,\ldots, 4$, and it is conjectured that no others exist; again, $f=0$ is excluded because the groups are solvable. In (B4), for $n$ to be prime it is necessary that $d$ should be prime; for $d=3$, over a hundred million such primes $n=p=(q^3-1)/(q-1)=q^2+q+1$, such as
\[2^2+2+1=7, \quad 3^2+3+1=13, \quad 5^2+5+1=31, \quad 8^2+8+1=73\]
and so on, are known to exist (see~\cite{JZ2}), and we conjecture, as would follow from a proof of the Bateman--Horn Conjecture (see Section~\ref{sec:NT}), that infinitely many such {\sl projective primes\/} $(q^d-1)/(q-1)$ exist for each prime $d\ge 3$. In (B5), for $n$ to be a prime power it is also necessary that $d$ should be prime (see~Corollary~\ref{cor:dprime}); at the time of writing, the only example we know where $n$ is a proper prime power is the group ${\rm PSL}_5(3)$, which has natural degree $11^2$; we conjecture, as would follow from a proof of the Nagell--Ljunggren Conjecture (see Section~\ref{sec:NT}), that there are no other examples. (The groups in (B4) and (B5) are listed separately, partly because of these very different conjectures, and partly because they are treated differently in our next theorem.) In (C) ${\rm PSL}_2(11)$ acts on the cosets of subgroups isomorphic to ${\rm A}_5$, the two conjugacy classes of these giving the two representations; ${\rm M}_{11}$ and ${\rm M}_{23}$ are Mathieu groups, automorphism groups of Steiner systems $S(4,5,11)$ and $S(4,7,23)$ with $11$ and $23$ points. In (D), $W(E_6)$ is the Weyl group (of the simple Lie group and Lie algebra) of type $E_6$, isomorphic to the group of the $27$ lines on a cubic surface, and $W(E_6)'$ is its commutator subgroup of index $2$, a simple group of order $25920$; they both have rank $3$, corresponding to the fact that pairs of lines on a cubic surface may be identical, intersecting or disjoint, giving three orbits. For more background on these groups, see~\cite{CCNPW} or \cite{Wil}.

There are two repeated representations in this list: the actions of ${\rm A}_5$ and ${\rm S}_5$ in (A1) can be identified with those of ${\rm PSL}_2(4)$ and ${\rm P\Sigma L}_2(4)={\rm PSL}_2(4)\rtimes{\rm C}_2$ in (B2). In addition, there are two instances of inequivalent representations of the same or isomorphic groups: ${\rm S}_{2^k}$ ($k\ge 3$) appears in (A1) in tts natural representation, and in (A2) as an imprimitive group of degree $2^{k+1}$, while ${\rm PSL}_2(7)$ appears in (B1) with degree $8$ and again in (B4) as the isomorphic group ${\rm PSL}_3(2)$ with degree $7$.

In this theorem, the only families of groups which are provably infinite are those in (A). It is conjectured that those in (B1) and (B4) are also infinite, and that those in (B2) and (B5) are finite; however, these conjectures involve number-theoretic problems which seem to be beyond the scope of current knowledge, and all one can do in these cases is present heuristic arguments and experimental evidence (see Section~\ref{sec:NT}).

A $p$-complement in a finite group $G$ is a subgroup of index a power of $p$ and order coprime to $p$; it is thus a sort of dual to a Sylow $p$-subgroup, where these conditions are reversed. Unlike Sylow subgroups, however, $p$-complements do not always exist, and when they do they need not all be conjugate, or even isomorphic (Buturlakin and Revin~\cite{BR} and Nesterov~\cite{Nes} have given examples of this; see Section~\ref{sec:conj} for details). If it exists, a $p$-complement is a point stabiliser in a transitive permutation representation of $G$ of $p$-power degree. Conversely, the point stabilisers in such a permutation group are $p$-complements provided they are $p'$-groups. It follows that the almost simple groups with $p$-complements form a subset of the groups in Theorem~\ref{th:mainppower}. They are classified in our second main result, which extends a result of Kazarin~\cite{Kaz} for simple groups:

\begin{theo}\label{th:mainpcomp}
The almost simple groups with $p$-complements, of index $n=p^k$, are as follows:
\begin{itemize}
\item[\rm(A)] ${\rm A}_n$ and ${\rm S}_n$ for $n=p\ge 5$;
\item[\rm(B)] 
\begin{enumerate}
\item ${\rm PSL}_2(q)$ and ${\rm PGL}_2(q)$ for any Mersenne prime $q=2^k-1>3$, with $p=2$;
\item ${\rm PSL}_2(q)\rtimes{\rm C}_{2^i}\le{\rm P\Gamma L}_2(q)$ for $q=2^e>2$ and $i=0,\ldots,f$, where $e=2^f$ and $p^k=p=2^e+1$ is a Fermat prime;
\item ${\rm PSL}_2(8)$ and ${\rm P\Gamma L}_2(8)$, where $p^k=9$;
\item for $d>2$ and $k=1$, any group $G$ such that ${\rm PSL}_d(q)\le G\le {\rm P\Gamma L}_d(q)$ where $n=(q^d-1)/(q-1)$ is a prime $p$;
\item for $d>2$ and $k>1$, any group $G$ such that ${\rm PSL}_d(q)\le G\le {\rm P\Gamma L}_d(q)$ where $n=(q^d-1)/(q-1)$ is a prime power $p^k$, and
$|G:G\cap{\rm PGL}_d(q)|$ is coprime to $p$;
\end{enumerate}
\item[\rm(C)] ${\rm PSL}_ 2(11)$ and ${\rm M}_{11}$ for $n=11$, and ${\rm M}_{23}$ for $n=23$.
\end{itemize}
There are two conjugacy classes of $p$-complements for $G$ in {\rm(B4)} and {\rm(B5)}, and ${\rm PSL}_2(11)$ in {\rm(C)}, and there is one in all the other groups.
In all cases, the $p$-complements are maximal subgroups.
\end{theo}

The groups in Theorem~\ref{th:mainppower} which are excluded here are all those in (A) for proper prime powers $n$, the two groups in (D), and possibly some groups in (B5) where $p$ divides
$|G:G\cap{\rm PGL}_d(q)|$: the only known group with $d>2$ and $k>1$, namely ${\rm PSL}_5(3)$, satisfies the additional necessary condition given here, and its point stabilisers are $p$-complements for $p=11$. Our earlier remarks about the cardinalities of the families of groups in Theorem~\ref{th:mainppower} also apply here.

In Section~\ref{sec:conj} we discuss the extension by Nesterov~\cite{Nes}, to a set $\mathcal NC$ of primes including the projective primes, of some examples due to Buturlakin and Revin~\cite{BR} for $p=7$. We strengthen his results by showing, for instance, that if $p\in{\mathcal NC}$ then a group can have an arbitrary number of mutually non-isomorphic $p$-complements.
After some preparatory work in Sections~\ref{sec:p-comp}--\ref{sec:permgps}, Theorems~\ref{th:mainppower} and \ref{th:mainpcomp} are proved, first for simple groups and then for almost simple groups, in Sections~\ref{sec:simple}--\ref{sec:almostB}.

The O'Nan--Scott Theorem (see~\cite[Chapter~4]{DM}) divides finite primitive permutation groups into a small number of classes, depending on their normal subgroup structure. By inspection, the only classes which can contain groups of prime power degree are those of almost simple, holomorph affine or product action type (see Section~\ref{sec:allprim} for an explanation of these types). This allows us to extend some of our results to primitive groups in general.

For example, there is current  interest (in connection with Cayley graphs and Schur rings, for instance) in the question of which permutation groups have regular subgroups. Those of holomorph affine type obviously do; by applying Theorem~\ref{th:mainppower} to the groups of almost simple and product action type (which are constructed from almost simple subgroups) we deduce the following:
\begin{coro}\label{cor:regsgp}
Every primitive permutation group of prime power degree has a regular subgroup.
\end{coro}

Theorem~\ref{th:mainppower} shows that if $G$ is an almost simple group then
\begin{itemize}
\item apart from $G\cong{\rm PSL}_3(2)\cong{\rm PSL}_2(7)$ of degrees $7$ and $8$, $G$ can have a faithful transitive permutation representation of $p$-power degree for at most one prime $p$;
\item $G$ has, up to equivalence, at most two such representations;
\item if there are two such representations then apart from $G={\rm S}_n$ ($n=2^k\ge 8$) of degrees $2^k$ and $2^{k+1}$ they have the same degree $p^k$ and are equivalent under an automorphism of $G$.
\end{itemize}
 It follows that a similar result applies to primitive groups of product action type. This also applies to primitive 
 holomorph affine groups $G=V\rtimes G_0\le{\rm AGL}_k(p)$, which have a regular normal translation subgroup $V$ and an irreducible point stabiliser $G_0\le{\rm GL}_k(p)$; however, in this case there is no upper bound on the number $|H^1(G_0,V)|$ of such representations of $G$, but nevertheless they are mutually equivalent under automorphisms of $G$. Thus we have the following:
\begin{coro}\label{cor:autequiv}
If $G$ is a finite group, then any two faithful primitive permutation representations of $G$ of the same prime power degree are equivalent under an automorphism of $G$.
\end{coro}

These two corollaries are proved in Sections~\ref{sec:regular} and \ref{sec:ineq}, where examples are given to show that they fail if the conditions of primitivity or prime power degree are relaxed. 

Any group with $p$-complements maps onto one with $p$-complements which are core-free, so that the corresponding action on cosets is faithful. Theorem~\ref{th:mainpcomp} shows which primitive permutation groups of almost simple type have $p$-complements; as with Theorem~\ref{th:mainppower} for regular subgroups, this result can also also applied to those of product action type, while the holomorphic affine groups are easily treated. This yields a classification of primitive groups with $p$-complements stated and proved in Section~\ref{sec:allprim}.

In the proofs of Theorems~\ref{th:mainppower} and \ref{th:mainpcomp}, cases (A), (C) and (D) are fairly straightforward; however case (B) involves a significant amount of number theory, some of it leading to difficult open problems and unproved conjectures in that area. These are discussed in Section~\ref{sec:NT}, placed at the end of the paper to avoid interrupting the flow of the group-theoretic arguments.


\section{$p$-complements}\label{sec:p-comp}

Throughout, $G$ is a finite group, and $p$ is a prime number dividing $|G|$.
A subgroup $H$ of $G$ is defined to be a $p$-{\sl complement\/} in $G$ if it has index a power of $p$ and order coprime to $p$,
that is, $H$ is a Hall $p'$-subgroup of $G$. In a sense, this is the dual of the concept of a Sylow $p$-subgroup, which has $p$-power order and $p'$ index.

Any group $H$ of order coprime to $p$ can appear as a $p$-complement: simply take $G$ to be a semidirect product $H\rtimes P$ of a normal subgroup $H$ by a $p$-group $P$. Conversely, by the Schur--Zassenhaus Theorem, every group with a normal $p$-complement has this form. In order to avoid trivial examples such as this, we will concentrate mainly on cases where the $p$-complement $H$ has trivial core in $G$, that is, it contains no non-identity normal subgroups of $G$.

By Hall's Theorems, every finite solvable group has a single conjugacy class of $p$-complements for each prime $p$ dividing its order (in fact, a single class of $\pi$-complements for each set $\pi$ of such primes).
 It is widely believed that as $x\to\infty$, almost all groups of order $n\le x$ are solvable (see~\cite{CEG} for evidence), 
so this would apply to them. Indeed, it is conjectured that almost all of them are $p$-groups for some prime $p$, and thus have the identity subgroup as a $p$-complement; for example, Higman~\cite{Hig} proved that there are at least $p^{2k^2(k-6)/27}$ groups of order $p^k$. 
In fact Erd\H os~\cite{Erd} proved that for almost all $n\le x$ every group of order $n$ has a normal Sylow $p$-subgroup for the largest prime $p$ dividing $n$, so that by the Schur--Zassenhaus Theorems it has a single conjugacy class of $p$-complements; specifically, he showed that the proportion of exceptional values of $n$ is at most 
\[\frac{1}{\exp\left(\bigl(\frac{1}{\sqrt 2}+O(1)\bigr)\bigl(\log x\log\log x\bigr)^{1/2}\right)}\to 0\quad\hbox{as}\quad x\to\infty.\]

Unlike Sylow's Theorems, these properties fail for non-solvable groups. For example, ${\rm A}_5$ has a single class of $5$-complements but no $2$- or $3$-complements, ${\rm PSL}_2(7)$ has one class of $2$-complements and two classes of $7$-complements, but no $3$-complements, while ${\rm A}_6$ has no $p$-complements for any prime~$p$. Indeed, our Theorem~\ref{th:mainpcomp} shows that this last phenomenon is typical of most finite almost simple groups.

For later use we include the following elementary results:

\begin{lemm}\label{lemma:KleG}
If $G$ has a subgroup $K$ of index a power of $p$, then any $p$-complement in $K$ is also a $p$-complement in $G$.
\end{lemm}

\begin{lemm}\label{lemma:KnornalG}
Let $K$ be a normal $p'$-subgroup of $G$, and let $K\le H\le G$. Then $H$ is a $p$-complement in $G$ if and only if $H/K$ is a $p$-complement in $G/K$.
\end{lemm}


\section{Conjugacy classes of $p$-complements}\label{sec:conj}

In~\cite{BR}, Buturlakin and Revin gave an example of a finite group with two non-isomorphic $7$-complements, and one in which the $7$-complements are isomorphic but are not conjugate in the automorphism group. These are based on the fact that the group ${\rm PSL}_3(2)$ has two conjugacy classes of $7$-complements, which are transposed by outer automorphisms. Their examples have been generalised by Nesterov~\cite{Nes} as follows. He defines, for each prime $p$, the following sets of groups:
\begin{itemize}
\item $\mathfrak{C}(p)$ is the set of finite groups $G$ such that all $p$-complements in $G$ are conjugate in $G$;
\item $\mathfrak{A}(p)$ is the set of finite groups $G$ such that all $p$-complements in $G$ are conjugate in ${\rm Aut}\,G$;
\item $\mathfrak{I}(p)$ is the set of finite groups $G$ such that all $p$-complements in $G$ are isomorphic;
\item $\mathfrak G$ is the set of all finite groups.
\end{itemize}
 Then clearly
\begin{equation}\label{eq:Nesterov}
\mathfrak{C}(p)\subseteq\mathfrak{A}(p)\subseteq\mathfrak{I}(p)\subseteq\mathfrak{G}.
\end{equation}

Nesterov now defines $\mathcal NC$ to be the set of primes $p$ such that
\begin{equation}\label{eq:NC}
p^k=\frac{q^d-1}{q-1}
\end{equation}
for some prime power $q$ and integers $d\ge 3$ and $k\ge 1$. For example, $2, 3, 5\not\in{\mathcal NC}$, but $7, 11, 13\in{\mathcal NC}$ since
\[7=\frac{2^3-1}{2-1},\quad 11^2=\frac{3^5-1}{3-1}\quad\hbox{and}\quad 13=\frac{3^3-1}{3-1}.\]
Taking $k=1$ shows that $\mathcal C$ contains all those projective primes $p$ for which $d\ge 3$. It also containes $p=11$ with $k=2$. In fact, $11$ is the only known element of $\mathcal NC$ with $k\ge 2$, and if the Nagell--Ljunggren Conjecture (see Section~\ref{sec:NT}) is true, it is unique.

Nesterov~\cite[Theorem 1]{Nes} has proved the following:

\begin{theo}[Nesterov]\label{th:Nesterov}
\begin{itemize}
\item[\rm(a)]If $p\in{\mathcal NC}$ then the inclusions in {\rm(\ref{eq:Nesterov})} are all proper.
\item[\rm(b)] If $p\not\in{\mathcal NC}$ then  the sets in {\rm(\ref{eq:Nesterov})} are all equal.
\end{itemize}
\end{theo}
 Here we will strengthen his statement (a) by showing that for each inclusion in~(\ref{eq:Nesterov}) there are groups in the larger set which fail to be in the proper subset by having, not just two, but an arbitrarily large number of inequivalent $p$-complements.

The significance of the integers $p^k$ appearing in (\ref{eq:NC}) is that these are the prime powers arising as the degrees of two transitive (in fact doubly transitive) permutation representations of $G={\rm PSL}_d(q)$ with $d\ge 3$, namely on points and and on hyperplanes of the associated projective geometry, giving two conjugacy classes of stabilisers which are $p$-complements in $G$ and are transposed in ${\rm Aut}\,G$. This immediately shows that the inclusion $\mathfrak{C}(p)\subseteq\mathfrak{A}(p)$ is proper for each $p\in{\mathcal NC}$. The following example extends this by showing that a group can have an unbounded number of conjugacy classes of $p$-complements, all equivalent under automorphisms.

\medskip

\noindent{\bf Example 3.1} If $p\in{\mathcal NC}$ then let $H={\rm PSL}_d(q)$ where $(q^d-1)/(q-1)=p^k$ and $d\ge 3$, and let $H_0$ and $H_1$ be the stabilisers in $H$ of a point and a hyperplane. For any $m\in{\mathbb N}$ let $G$ be the cartesian product $H^m$ of $m$ copies of $H$. Then $G$ has $2^m$ conjugacy classes of $p$-complements, represented by the subgroups $H_{\phi}=H_{\phi(1)}\times\cdots\times H_{\phi(m)}$ for each function $\phi:\{1,\ldots,m\}\to\{0,1\}$. Since $H_0$ and $H_1$ are conjugate in ${\rm Aut}\,H$, it follows that these $p$-complements $H_{\phi}$ are all conjugate in ${\rm Aut}\,G$.

\medskip

Similarly, in the case of the proper inclusion $\mathfrak{A}(p)\subset\mathfrak{I}(p)$, the following example shows that a group can have all of its $p$-complements mutually isomorphic, but lying in an unbounded number of orbits of its automorphism group.

\medskip

\noindent{\bf Example 3.2} With the notation of Example~3.1, let $G=H\wr{\rm C}_m=H^m\rtimes{\rm C}_m$ where now $m=p^e$ for some integer $e\ge 1$; the complement ${\rm C}_m$ acts by conjugation on the base group $H^m$ by cyclically shifting coordinates. The $2^m$ subgroups $H_{\phi}$ are mutually isomorphic $p$-complements in $G$. Such subgroups $H_{\phi}$ and $H_{\psi}$ are conjugate in $G$ if and only if $\phi$ and $\psi$ differ by a cyclic permutation of their domain. Since $m=p^e$, each function $\phi$ has $p^f$ images under such transformations for some $f=0, 1, \ldots, e$, so $H_{\phi}$ is in a conjugacy class of size $p^{m+f}$ in $G$. All such values of $f$ are realised by some $\phi$, so since automorphisms permute conjugacy classes of subgroups and preserve their size it follows that ${\rm Aut}\,G$ has at least $e+1$ orbits on the $p$-complements $H_{\phi}$. (In fact, in these cases there are many more than $e+1$ orbits: the exact number, which is asymptotic to $2^{m-1}/m$ as $m\to\infty$, can be obtained from counting two-coloured necklaces of $m$ beads.)

\medskip

In the case of the proper inclusion $\mathfrak{I}(p)\subset\mathfrak{G}$, the next example, generalising examples in~\cite{BR} and~\cite{Nes}, shows that a group can have an arbitrary number of isomorphism classes of $p$-complements.

\medskip

\noindent{\bf Example 3.3} First we construct a group with two non-isomorphic $p$-complements. Continuing the earlier notation, but lifting back from projective to linear groups, let $\tilde H={\rm GL}_d(q)$, and for $i=0$ or $1$ let $\tilde H_i$ be the inverse image of $H_i$ in $\tilde H$, that is, the subgroup stabilising a subspace $U_i$ of dimension $1$ or $d-1$ in the vector space $V={\mathbb F}_q^d$. Since $|\tilde H_i|=(q-1)|H_i|$ with $q-1$ coprime to $p$, these two subgroups $\tilde H_i$ are $p$-complements in $\tilde H$; they are conjugate in ${\rm Aut}\,\tilde H$ but not in $\tilde H$.

Now let $\overline H_i$ be the subgroup $\langle V, \tilde H_i\rangle=V\rtimes \tilde H_i$ of $\overline H:={\rm AGL}_d(q)=V\rtimes \tilde H$, for $i=0, 1$. Since $|V|=q^d$ is coprime to $p$, each $\overline H_i$ is a $p$-complement in $\overline H$. However, $\overline H_0\not\cong \overline H_1$. One can see this by checking that the only non-trivial normal subgroup of $\overline H_i$ properly contained in $V$ is $U_i$, of order $q$ or $q^{d-1}$ as $i=0$ or $1$. Thus $\overline H$ has two non-isomorphic $p$-complements (the example in~\cite{BR} has $\overline H={\rm AGL}_3(2)$ and $p=7$).

Finally, for any $m\in{\mathbb N}$ let $G$ be the cartesian product $\overline H^m$ of $m$ copies of $\overline H$. The functions $\phi:\{1,\ldots, m\}\to\{0, 1\}$ determine $2^m$ mutually non-conjugate $p$-complements $\overline H_{\phi}=\overline H_{\phi(1)}\times\cdots\times \overline H_{\phi(m)}$ in $G$. We have $\overline H_{\phi}\cong \overline H_{\psi}$ if and only if $|\phi^{-1}(1)|=|\psi^{-1}(1)|$, and this number can take any value from $0$ to $m$, so $G$ has $m+1$ isomorphism classes of $p$-complements. 

\medskip

In order to prove statement (b) of Nesterov's Theorem~\ref{th:Nesterov}, one needs to show that if $p\not\in{\mathcal NC}$ then in any finite group with $p$-complements, they form a single conjugacy class. Using results of Chunikhin, an argument by inductiion on the length of a composition series reduces the problem to the case of a nonabelian finite simple group, and here one can verify the claim by inspecting Guralnick's list~\cite{Gur} of those with transitive permutation representations of $p$-power degree. The fact that $p\not\in{\mathcal NC}$ means that the exceptions ${\rm PSL}_d(q)$ with $d\ge 3$, or with $d=2$ and $q=11$, are avoided.


\section{$p$-complements and permutation groups}\label{sec:p-compperm}

The following result gives alternative characterisations of $p$-complements:

\begin{prop}\label{prop:permgps}
If $H$ is a subgroup of $G$, and $P$ is a Sylow $p$-subgroup of $G$, then the following are equivalent:
\begin{itemize}
\item $H$ is a $p$-complement in $G$;
\item $H$ is a point stabiliser in a transitive permutation representation of $G$ of $p$-power degree, and $H$ is a $p'$-group;
\item $G=PH$ and $P\cap H=1$;
\item $P$ acts regularly on the cosets of $H$ in $G$;
\item $H$ acts regularly on the cosets of $P$ in $G$.
\end{itemize}
\end{prop}

The proof is elementary, using the following:
\begin{lemm}
If $G$ acts transitively on a set $\Omega$ with $|\Omega|=p^k$ then each Sylow $p$-subgroup $P$ of $G$ is also transitive on $\Omega$.
\end{lemm}

\noindent{\sl Proof.} Let $\alpha\in\Omega$. Then $|G:G_{\alpha}|\cdot|G_{\alpha}:P_{\alpha}|=|G:P_{\alpha}|=|G:P|\cdot|P:P_{\alpha}|$ with $|G:G_{\alpha}|=p^k$ and $|G:P|$ coprime to $p$, so $|P:P_{\alpha}|$ is divisible by $p^k$. Thus $\alpha$ lies in a $p$-orbit of size at least $p^k$, which must be $\Omega$. \hfill$\square$

\medskip

By Proposition~\ref{prop:permgps}, a group $G$ has $p$-complements if and only if there is a normal $p'$-subgroup $K$ of $G$ such that $G/K$ is a transitive permutation group of $p$-power degree with $p'$-groups as point stabilisers, so in order to study $p$-complements it is useful to consider permutation groups of prime power degree. For instance, the following result gives some examples:

\begin{coro}\label{cor:primedegree}
If $G$ is a transitive permutation group of prime degree $p$ then a point stabiliser $H$ is a $p$-complement in $G$.
\end{coro}

\noindent{\sl Proof.} One can identify $G$ and $H$ with subgroups of the symmetric groups ${\rm S}_p$ and ${\rm S}_{p-1}$. The latter has order $(p-1)!$, so $H$ is a $p'$-group, of index $p$ in $G$. \hfill$\square$

\medskip

As shown by Galois the solvable transitive groups of prime degree $p$ are the subgroups of the affine group ${\rm AGL}_1(p)$ containing the translation subgroup, one for each divisor of $p-1$. It follows from a result of Burnside~\cite{Bur95} and the classification of finite simple groups (see~\cite{Cam}, for example), that the nonsolvable groups are as follows: 
\begin{prop}\label{prop:degp}
The nonsolvable transitive permutation groups of prime degree $p$ are the following:
\begin{itemize}
\item[\rm(A)] ${\rm A}_p$ and ${\rm S}_p$ where $p\ge 5$,
\item[\rm(B)] groups $G$ such that ${\rm PSL}_d(q)\le G\le{\rm P\Gamma L}_d(q)$ where the natural degree $(q^d-1)/(q-1)$ is a prime $p$, and
\item[\rm(C)] ${\rm PSL}_2(11)$ and ${\rm M}_{11}$ where $p=11$, and ${\rm M}_{23}$ where $p=23$.
\end{itemize}
\end{prop}

See the Introduction for an explanation of these groups and their actions. It is an open problem whether there are finitely or infinitely many groups $G$ in (B): see~\cite{JZ2} and Section~\ref{sec:NT}.

A $p$-complement $H$ is a maximal subgroup of a group $G$  if and only if the corresponding permutation group is primitive. If it is not, then subgroups $L$ with $H<L<G$ have smaller $p$-power index in $G$, and correspond bijectively to systems of imprimitivity for $G$, which have $|G:L|$ blocks of size $|L:H|$.


'\section{Permutation groups of prime power degree}\label{sec:permgps}

Any transitive group of $p$-power degree maps onto a primitive group of (possibly smaller) $p$-power degree, so it is useful to concentrate on the latter.
The O'Nan--Scott Theorem (see~\cite[Chapter~4]{DM}) divides finite primitive permutation groups into five types: holomorph affine, almost simple, product action, simple diagonal and twisted wreath. Groups in the last two classes have degree a power of the order of a non-abelian finite simple group, and this cannot be a prime power. Primitive holomorph affine groups $G$ have the form $V\rtimes H$, where $V$ is a $k$-dimensional vector space over ${\mathbb F}_p$, with an irreducible point stabiliser $H\le{\rm GL}_k(p)$; this is a $p$-complement if and only if it is a $p'$-group, in which case the Schur--Zassenhaus Theorem implies that there is just one conjugacy class of $p$-complements.
Primitive groups of product action type and of $p$-power degree have socle $T^d$ where $T$ is a non-abelian finite simple group of smaller $p$-power degree, so this reduces many questions about these groups to questions about almost simple groups (see Section~\ref{sec:allprim} for further details). We will therefore concentrate on almost simple groups of prime power degree for the first part of this paper.

Any transitive permutation group $G$ of $p$-power degree either acts primitively, or has a primitive representation of smaller $p$-power degree; if $G$ is almost simple, this primitive representation is faithful. This enables us to extend our study of almost simple groups of prime power degree from primitive groups to transitive groups.
From now on we will assume that $G$ is an almost simple transitive permutation group of prime power degree $n=p^k$.


\section{Simple groups}\label{sec:simple}

We will first take $G=S$, a non-abelian finite simple group. The following result is due to Guralnick~\cite{Gur} (see also~\cite[Theorem~5.8]{AF}):

\begin{theo}\label{th:Guralnick}
A non-abelian finite simple $S$ has a transitive permutation representation of prime power degree $n=p^k$ if and only if
\begin{itemize}
\item[\rm(A)]  $S={\rm A}_n$ where $n=p^k\ge 5$, or
\item[\rm(B)] $S={\rm PSL}_d(q)$ where $n=(q^d-1)/(q-1)=p^k$, or
\item[\rm(C)] $S={\rm PSL}_2(11)$ or ${\rm M}_{11}$ where $n=11$, or $S={\rm M}_{23}$ where $n=23$, or
\item[\rm(D)] $S={\rm PSp}_4(3)=W(E_6)'$ where $n=27$.
\end{itemize}
\end{theo}

Again, see the Introduction for an explanation of these groups and their actions. In (B), ${\rm PSL}_d(q)$ is simple if and only if either $d\ge 3$, or $d=2$ and $q\ge 4$. It is a difficult number-theoretic problem to determine the set of pairs $d$ and $q$ for which the degree $n=(q^d-1)/(q-1)$ is a prime power, although of course, it is straightforward in individual cases: indeed, it is not known whether this set is finite or infinite (see Section~\ref{sec:NT} for more on this problem). One necessary condition is that $d$ should be prime: when $k=1$ an obvious factorisation of $q^{d-1}+q^{d-2}+\cdots+q+1$ shows that $d$ cannot be composite, whereas the proof for $k>1$, given in Corollary~\ref{cor:dprime}, relies on Zsigmondy's Theorem.

In this list, a point stabiliser $S_{\alpha}$ is $p$-complement if and only if it is a $p'$-group, or equivalently $|S|$ is not divisible by $p^{k+1}$. Clearly, this happens in (A) if and only if $n$ is prime, in all three groups in (C), but not in (D) because here $p=3$ and $|S_{\alpha}|=960$. The situation in (B) is more complicated, and will be revisited in Section~\ref{sec:almostB}, where we will show that all such groups have $p$-complements. To summarise, we have the following:

\begin{coro}\label{cor:simplepcomp}
A non-abelian finite simple group $S$ has a $p$-complement if and only if
\begin{itemize}
\item[\rm(A)]  $S={\rm A}_n$ where $n=p\ge 5$, or
\item[\rm(B)] $S={\rm PSL}_d(q)$ where $(q^d-1)/(q-1)=p^k$ and $p^{k+1}$ does not divide $|S|$, or
\item[\rm(C)] $S={\rm PSL}_2(11)$ or ${\rm M}_{11}$ where $p=11$, or $S={\rm M}_{23}$ where $p=23$.
\end{itemize}
There are two conjugacy classes of $p$-complements for ${\rm PSL}_2(11)$ in {\rm(C)}, and two for ${\rm PSL}_d(q)$ in {\rm(B)} if $d\ge 3$; these classes are transposed in ${\rm Aut}\,S$ in both cases. All the other groups $S$ have just one class of $p$-complements. 
\end{coro}

In fact, we will show in Propositions~\ref{prop:d=2} and~\ref{prop:PSLp-comp} that the final condition in (B) is redundant. Note that ${\rm A}_5\cong {\rm PSL}_2(4)$, with $p=5$ in each case. We also have ${\rm PSL}_3(2)\cong{\rm PSL}_2(7)$, with $p=7$ and $p=2$ respectively; there are no other isomorphisms between these groups, so this is the only example of a simple group with $p$-complements for two different primes $p$.

We will restrict our attention now to case~(B), the only case presenting any difficulties, namely, the questions of whether the natural degree
\begin{equation}\label{eq:degree}
n=\frac{q^d-1}{q-1}=q^{d-1}+q^{d-2}+\cdots+q+1
\end{equation}
of $S$ is a prime power $p^k$, and if so whether the point stabilisers are $p'$-groups, the two conditions in Corollary~\ref{cor:simplepcomp}(B).

We first consider the case $d=2$, so that $S={\rm PSL}_2(q)$ with $n=q+1$ and we need to know whether this has the form $p^k$ for some prime $p$. First we need a simple lemma.

\begin{lemm}\label{lemma:8and9}
If $p^k=r^1+1$, where $p$ and $r$ are prime numbers and $k, l>1$, then $p^k=9$ and $r^l=8$.
\end{lemm}

\noindent{\sl Proof.} Since $p^k$ and $r^l$ have opposite parities, either $p$ or $q$ is $2$, and the other is odd.
If $p=2$ then $r^l=2^k-1\equiv-1$ mod~$(4)$, so $r\equiv -1$ mod~$(4)$ and $l$ is odd. Then $2^k=r^l+1=(r+1)(r^{l-1}-r^{l-2}+\cdots-r+1)$ with the second factor having an odd number of odd terms, so that it is odd; this is impossible since it divides $2^k$ and is greater than $1$.

If $r=2$ then $p^k=2^l+1\equiv 1$ mod~$(4)$, so either $k$ is even, or $k$ is odd and $p\equiv 1$ mod~$(4)$. In the first case, writing $k=2m$ gives $2^l=(p^m+1)(p^m-1)$, so $p^m\pm 1$ are both powers of $2$ and hence $p^m=3$, $p^k=9$ and $r^l=8$ as required. In the second case $2^l=p^k-1=(p-1)(p^{k-1}+p^{k-2}+\cdots+p+1)$ with the second factor odd, leading to a contradiction as before. \hfill$\square$

\medskip

(Of course, one can deduce this lemma from Mih\u ailescu's proof~\cite{Mih} of the Catalan Conjecture, which asserts the same result but without the hypothesis that $p$ and $r$ are prime; however, such a deep result is not needed here.)

\begin{prop}\label{prop:d=2}
The group $S={\rm PSL}_2(q)$, acting naturally, has prime power degree $p^k$ if and only if either
\begin{enumerate}
\item $p=2$ and $q=2^k-1$ is a Mersenne prime, or
\item $k=1$ and $p=q+1=2^e+1$ is a Fermat prime, or
\item $p^k=9$ and $q=8$.
\end{enumerate}
In each case $S$ has a single conjugacy class of $p$-complements, namely the natural point stabilisers; they have order $q(q-1)/2$ in case~(1), and order $q(q-1)$ in cases~(2) and (3).
\end{prop}

\noindent{\sl Proof.} The natural degree of $S$ is $q+1$, so suppose that $q+1=p^k$. By Lemma~\ref{lemma:8and9} either $q$ or $p^k$ is prime, or these are respectively $8$ and $9$. If $q$ is prime, either $q=2$ and hence $p^k=p=3$, an instance of conclusion~(2), or $q$ is a Mersenne prime $p^k-1$ with $p=2$, as in conclusion~(1). If $p^k$ is prime then $k=1$ and $p=q+1$ is a Fermat prime, as in conclusion~(2). Finally, if $q=8$ and $p^k=9$ we have conclusion~(3).
Using the fact that $|S|=q(q^2-1)$ or $q(q^2-1)/2$ as $q$ is even or odd, we see that the point stabilisers are $p'$-groups of the stated order, so they are $p$-complements in $S$.
\hfill$\square$

\medskip

Note that ${\rm PSL}_2(4)\cong{\rm PSL}_2(5)$, so this group has two natural actions, of degrees $5$ and $6$; the first is isomorphic to its action as ${\rm A}_5$, whereas the second, not having prime power degree, is irrelevant here.
There seems to be no hope of classifying the Mersenne or Fermat primes (see the comments in the Introduction), so there is nothing more we can add about the case $d=2$.

We will now consider the case $k=1$, so that the degree $n$ is a prime $p$. Then $d$ must be prime, otherwise there is an obvious factorisation in (\ref{eq:degree}). We have dealt with the case $d=2$. The Bateman--Horn Conjecture~\cite{BH} suggests that for each prime $d\ge 3$ there are infinitely many {\sl projective primes\/} $p$ of this form, with prime powers $q$ (and even with prime $q$). More precisely, this conjecture gives formulae for the asymptotic distribution of projective primes, and these agree remarkably closely with what is found by computer search and primality testing; if true, the conjecture would imply the existence of infinitely many such primes for each prime $d\ge 3$. As an example, there are $129\,294\,308$ primes $q\le 10^{11}$ such that $q^2+q+1$ is prime. See Section~\ref{sec:NT} and~\cite{JZ2} for more details. 

We now consider whether the degree $n=p^k$ in Case~(B) can be a proper prime power, that is, $k>1$. If $d=2$ then the only possibilities are cases (1) and (3) of Proposition~\ref{prop:d=2}, all with one class of $p$-complements. The only known example with $d\ge 3$ is $n=11^2=(3^5-1)/(3-1)$ for ${\rm PSL}_5(3)$, with two classes of 11-complements, the stabilisers of points and hyperplanes.

The Bateman--Horn Conjecture, which concerns prime values of polynomials, does not apply here, but another number-theoretic conjecture is applicable.
It is conjectured that the only solutions $x, y\in{\mathbb Z}$ of the Nagell--Ljunggren equation
\begin{equation}\label{eq:NL}
\frac{x^d-1}{x-1}=y^k
\end{equation} 
with $|x|$, $|y|$, $k>1$ and $d>2$ are
\begin{equation}\label{eq:NLknown}
\frac{3^5-1}{3-1}=11^2,\quad \frac{7^4-1}{7-1}=20^2,\quad \frac{18^3-1}{18-1}=7^3\quad\hbox{and}\quad \frac{(-19)^3-1}{-19-1}=7^3.
\end{equation}
(See the survey~\cite{BM} by Bugeaud and Mignotte, and also~\cite[A208242]{OEIS}.) This conjecture is open, but it has been proved in some special cases (see Section~\ref{sec:NT}). Only the first solution, corresponding to $S={\rm PSL}_5(3)$, has $x$ and $y$ positive prime powers, so the other three are irrelevant here. We therefore conjecture that the only group ${\rm PSL}_d(q)$ of proper prime power degree for $d\ge 3$ is ${\rm PSL}_5(3)$.

The following is a special case, appropriate for applications to finite groups, of a more general theorem of Zsigmondy~\cite{Zsi}:

\begin{prop}[Zsigmondy's Theorem]\label{prop:Zsigmondy}
For any integers $q, d\ge 2$, one of the following holds:
\begin{itemize}
\item[{\rm(a)}] there is a prime $s$ dividing $q^d-1$ but not dividing $q^e-1$ for any $e<d$;
\item[{\rm(b)}] $q=2$ and $d=6$;
\item[{\rm(c)}] $q$ is a Mersenne prime and $d=2$.
\end{itemize}
\end{prop}

A prime $s$ satisfying (a) is called a {\sl primitive divisor} of $q^d-1$. Statement (a) fails in case~(b), since $2^6-1=3^2\cdot 7$ with $3$ dividing $2^2-1$ and $7$ dividing $2^3-1$, and it fails in case~(c) since if $q=2^a-1$ then $q^2-1=(q+1)(q-1)=2^a(q-1)$, so any prime dividing $q^2-1$ also divides $q-1$.

\begin{coro}\label{cor:dprime}
Let $q, d\ge 2$ and let $n=(q^d-1)/(q-1)$ be a prime power $p^k$. Then
\begin{itemize}
\item[{\rm(a)}] $d$ is prime;
\item[{\rm(b)}] if $d>2$ then $d$ divides $p-1$.
\end{itemize}
\end{coro}

\noindent{\sl Proof.} We may assume that $d>2$. Then Zsigmondy's Theorem implies that either $q=2$ and $d=6$, or some prime $s$ divides $q^d-1$ but not $q^e-1$ for any $e<d$. We can ignore  the first case since then $n=63$, which is not a prime power.
In the second case, since $s$ divides $q^d-1$ but not $q-1$, it divides $n$, so $s=p$; now if $d=ef$ with $e, f<d$ then $(q^e-1)/(q-1)$ divides $n$, so it is a power of $p$, and hence $s$ divides $q^e-1$, a contradiction. Thus $d$ is prime. Since $p$ divides $q^d-1$ but not $q-1$, this prime $d$ is the multiplicative order of $q$ mod~$(p)$, so it divides $|{\mathbb Z}_p^*|=p-1$. \hfill$\square$

\medskip

Note that Corollary~\ref{cor:dprime}(b) fails if $d=2$, by case~(1) of Proposition~\ref{prop:d=2}, where $p=2$. In fact, Corollary~\ref{cor:dprime}(b) and Proposition~\ref{prop:d=2} immediately imply the following:

\begin{coro}\label{cor:p=2}
If $p=2$ then $d=2$, and we have case~(1) of Proposition~\ref{prop:d=2}.
\end{coro}

\medskip

We may therefore assume from now on that $d$ and $p$ are odd primes.


{\section{$p$-complements in ${\rm PSL}_d(q)$}\label{sec:p-compPSL}

We have already considered the problem of determining when the natural degree $n$ of the group $S={\rm PSL}_d(q)$ is a prime power~$p^k$. If it is, then a point-stabiliser $S_{\alpha}$ is a $p$-complement in $S$ if and only if it is a $p'$-group. We have seen that the point stabilisers are $p$-complements when $d=2$, so we may assume that $d\ge 3$.

The group $S$ has order
\begin{equation}\label{eq:PSL}
|S|=\frac{1}{(d,q-1)}\cdot q^{d(d-1)/2}\cdot (q-1)^{d-1}\cdot\prod_{i=1}^{d-1}(q^i+q^{i-1}+\cdots+q+1).
\end{equation}
Its natural degree is
\[n=q^{d-1}+q^{d-2}+\cdots+q+1,\]
the factor of greatest degree in the final product in equation~(\ref{eq:PSL}), so $S_{\alpha}$ has order
\begin{equation}\label{eq:stab}
|S_{\alpha}|=\frac{|S|}{n}=\frac{1}{(d,q-1)}\cdot q^{d(d-1)/2}(q-1)^{d-1}\cdot\prod_{i=1}^{d-2}(q^i+q^{i-1}+\cdots+q+1).
\end{equation}
The problem therefore is to determine whether or not this is divisible by $p$ in those cases when the degree $n$ is a power of $p$.

\medskip

\noindent{\bf Example 7.1} If $S={\rm PSL}_5(3)$ then 
\[|S|=\frac{1}{(5,2)}\cdot 3^{10}\cdot 2^4\cdot\prod_{i=1}^4(3^i+3^{i-1}+\cdots+3+1)=3^{10}\cdot 2^4\cdot 121\cdot 40\cdot 13\cdot 4=2^9\cdot3^{10}\cdot5\cdot11^2\cdot13,\]
so $S_{\alpha}$ has order $2^9\cdot3^{10}\cdot5\cdot13$ and is therefore an $11$-complement, belonging to one of two conjugacy classes.

\medskip

\begin{prop}\label{prop:PSLp-comp}
If $S={\rm PSL}_d(q)$ has natural degree $n=p^k$ for some prime $p$ then the point stabilisers are $p$-complements.
\end{prop}

\noindent{\sl Proof.} We may assume that $d>2$. We need to show that $p$ divides none of the factors of $|S_{\alpha}|$ on the right-hand side in (\ref{eq:stab}). We have $q=r^e$ for some prime $r$, and we will apply Zsigmondy's Theorem to $q^d-1=r^{de}-1$. Of the two Zsigmondy exceptional cases, $de=2$ (and $r$ Mersenne) cannot arise, since we are assuming that $d>2$; in the other exceptional case, if $r^{de}=2^6$ then either $d=3$ and $q=4$, or $d=6$ and $q=2$, giving $n=21$ or $n=63$, neither of which is a prime power. It follows that there is a primitive divisor of $r^{de}-1$; this is a prime which divides $q^d-1$ but not $q-1$, so it divides $n$, and hence must be $p$. It cannot divide factors of $|S_{\alpha}|$ of the form $q-1=r^e-1$ or $q^i+\cdots+q+1=(q^i-1)/(q-1)$ for $i\le d-2$, and it cannot divide $q$ since $p^k=n\equiv 1$ mod~$(q)$. Thus $S_{\alpha}$ is a $p'$-group and therefore a $p$-complement. \hfill$\square$

\medskip

The following result, which follows immediately from Corollary~\ref{cor:dprime}(b), will be useful when we consider almost simple groups with socle ${\rm PSL}_d(q)$.

\begin{coro}\label{cor:d,q-1}
If $d>2$ and $S={\rm PSL}_d(q)$ has natural degree $n=p^k$ for some prime $p$ then $p$ does not divide $(d,q-1)$.
\end{coro}


\section{Almost simple groups}
\label{sec:almost}

Recall that $G$ is almost simple if $S\le G\le A:={\rm Aut}\,S$ for some non-abelian finite simple group $S$.
As in the case of simple groups $S$, we are interested in whether $G$ is a permutation group of $p$-power degree for some prime $p$, and if so whether it has $p$-complements, and how many conjugacy classes they form.

If $G$ acts transitively with degree $p^l$ on a set $\Omega$, its socle $S$ is transitive of degree $p^k$ for some $k\le l$: here either $k=l$, so that $S$ is transitive on $\Omega$, or $1\le k<l$, so that $S$ is intransitive on $\Omega$, with orbits of length $p^k$ forming blocks of imprimitivity for $G$, and with $p^{l-k}$ dividing $|G:S|$. In either case, $S$ is one of the groups of degree $p^k$ listed earlier in Theorem~\ref{th:Guralnick}. The subgroups $G$ of $A$ containing $S$ correspond bijectively to the subgroups $G/S$ of the outer automorphism group $A/S={\rm Out}\,S$, and this group is known for all finite simple groups $S$ (see \cite{CCNPW} or \cite{Wil}), so for each group $S$ in Theorem~\ref{th:Guralnick} one can check all such groups $G$ for having the relevant properties. Here we will do this for cases (A), (C) and (D) in turn, leaving the more difficult case (B) until the next section.

\medskip

(A) Here $S={\rm A}_n$ where $n=p^k\ge 5$. Since $n\ne 6$ we have $A={\rm S}_n$, with $|{\rm S}_n:{\rm A}_n|=2$, so only $G={\rm S}_n$ arises. One possibility is that $l=k$ and ${\rm S}_n$ acts naturally with degree $n$ on the cosets of ${\rm S}_{n-1}$, extending the action of ${\rm A}_n$ on cosets of ${\rm A}_{n-1}$, as in Theorem~\ref{th:mainppower}(A1); then ${\rm S}_n$ has $p$-complements if and only if $k=1$, in which case they form a single conjugacy class. However, if $p=2$ there is also the possibility that $l=k+1$, so that ${\rm S}_n$ acts with degree $2^l=2n$ on the cosets of ${\rm A}_{n-1}$, as in Theorem~\ref{th:mainppower}(A2), and ${\rm A}_n$ has two natural orbits in this action; now ${\rm S}_n$ has no $2$-complements, since ${\rm A}_{n-1}$ has even order.

\medskip

(C) Since ${\rm Out}\,S=1$ for the Mathieu groups $S={\rm M}_{11}$ and ${\rm M}_{23}$, no new groups $G$ arise from them. The group $S={\rm PSL}_2(11)$ has index $2$ in its automorphism group $A={\rm PGL}_2(11)$, so only $G=A$ arises; this has no subgroups of index~$11$ (see~\cite{CCNPW}), so it does not act transitively with degree a power of $11$.

\medskip

(D) Similarly $S={\rm PSp}_4(3)=W(E_6)'$ has index $2$ in $A=W(E_6)$, so only $G=W(E_6)$ arises. This is the full group of the $27$ lines, extending the action of $S$, so it also has degree $3^3$. However, as with $S$ the point stabilisers are not $3'$-groups, so it has no $3$-complements. Since $|G:S|=2$ is not divisible by $p=3$ there is no possibility of an intransitive action of $S$ extending to a transitive action of $G$ of degree a power of $3$.

\medskip

To summarise, we have the following:

\begin{prop}\label{prop:transnotB}
Apart from groups with socle $S={\rm PSL}_d(q)$, the almost simple groups with transitive representations of prime power degree $n=p^k$ are as follows:

\smallskip

{\rm (A)} ${\rm A}_n$ and ${\rm S}_n$ where $n=p^k\ge 5$, and ${\rm S}_{n/2}$ with $n=2^k\ge 16$;

\smallskip

{\rm (C)} ${\rm PSL}_2(11)$ and ${\rm M}_{11}$ where $n=11$, and ${\rm M}_{23}$ where $n=23$;

\smallskip

{\rm (D)} $W(E_6)$ and $W(E_6)'$ where $n=27$.

\smallskip

\noindent Each group has a single representation, apart from ${\rm PSL}_2(11)$ which has two.
Those groups with $p$-complements are as follows:

\smallskip

{\rm (A)} ${\rm A}_n$ and ${\rm S}_n$ where $n=p\ge 5$;

\smallskip

{\rm (C)} ${\rm PSL}_2(11)$ and ${\rm M}_{11}$ where $p=11$, and ${\rm M}_{23}$ where $p=23$.
\smallskip

\noindent Each group has a single conjugacy class of $p$-complements, apart from ${\rm PSL}_2(11)$ which has two.

\end{prop}

Almost simple groups with socle $S={\rm PSL}_d(q)$ are dealt with in the next section.


\section{Almost simple groups in case (B)}\label{sec:almostB}

In case~(B) we have $S={\rm PSL}_d(q)$, acting with degree $n=(q^d-1)/(q-1)$ on the points of ${\mathbb P}^{d-1}({\mathbb F}_q)$, and also on the hyperplanes if $d>2$. If $d=2$ then $A={\rm P\Gamma L}_d(q)$, extending this action of $S$, whereas if $d>2$ then $A={\rm P\Gamma L}_d(q)\rtimes{\rm C}_2$ with the complement ${\rm C}_2$ generated by the graph automorphism, transposing points and hyperplanes. We are interested in almost simple groups $G$ such that $S<G\le A$. If $d=2$ then every such group $G$ extends the natural action of $S$, but if $d>2$ then $G$ extends the actions of $S$ if and only if $G\le {\rm P\Gamma L}_d(q)$. If $G\not\le {\rm P\Gamma L}_d(q)$ with $d>2$ then $G$ acts on the set of points and hyperplanes with degree $2n$, and this is not a $p$-power unless $p=2$; but then Corollary~\ref{cor:p=2} implies that $d=2$, against our assumption. Hence, for each $d\ge 2$, we may restrict our attention to subgroups $G\le{\rm P\Gamma L}_d(q)$. Since such groups extend the action of $S$, they have prime power degree if and only if $S$ does, so this, together with Proposition~\ref{prop:transnotB}, completes the proof of Theorem~\ref{th:mainppower}.

In order to complete the proof of Theorem~\ref{th:mainpcomp} we need to consider which subgroups $G\le{\rm P\Gamma L}_d(q)$, with socle $S$ of prime power degree $p^k$, have point stabilisers which are $p'$-groups and are thus $p$-complements. When $d=2$ we have a solution:

\begin{prop}\label{prop:ASd=2}
If $S={\rm PSL}_2(q)$ has degree $p^k$ for some prime $p$ then every almost simple group $G$ with socle $S$ has a single conjugacy class of $p$-complements.
In the various cases of Proposition~\ref{prop:d=2} the groups $G$ (excluding $S$) are:
\begin{enumerate}
\item ${\rm PGL}_2(q)={\rm PSL}_2(q)\rtimes{\rm C}_2$, where $q=2^k-1$ is a Mersenne prime and $p=2$;
\item ${\rm PSL}_2(q)\rtimes{\rm C}_{2^i}\le {\rm P\Gamma L}_2(q)$ for $i=1,\ldots,f$ where $q=2^e$, $e=2^f$ and $p^k=p=2^e+1$ is a Fermat prime;
\item ${\rm P\Gamma L}_2(8)={\rm PSL}_2(8)\rtimes{\rm C}_3$, where $p^k=9$.
\end{enumerate}
\end{prop}

\noindent{\sl Proof.} Since $d=2$, $S$ is as in case~(1), (2) or (3) of Proposition~\ref{prop:d=2}. In each of these cases $S$ has a single conjugacy class of $p$-complements, of order $q(q-1)/2$ in case~(1), and $q(q-1)$ in cases~(2) and (3). In case~(1) only $G=A={\rm PGL}_2(q)$ arises; since $p=2=|G:S|$ the $2$-complements in $S$ are also $2$-complements in $G$ by Lemma~\ref{lemma:KleG}. Case~(3) is similar, but with $p=3=|G:S|$. In case~(2), each $G\le A$ is a transitive group of prime degree $p$, so its point stabilisers are $p$-complements by Corollary~\ref{cor:primedegree}. \hfill$\square$

\medskip

Now suppose that $S\le G\le{\rm P\Gamma L}_d(q)$ where $d>2$. If $k=1$ then $G$ is a transitive group of prime  degree $p$, so it has $p$-complements by Corollary~\ref{cor:primedegree}, and these form two conjugacy classes, the stabilisers of points and of hyperplanes. Hence we may assume that $k>1$.

\medskip

\noindent{\bf Example 9.1} In the only known example where $d>2$ and $n$ is a proper prime power, namely $S={\rm PSL}_5(3)$ with $n=11^2$, the stabilisers of points and of hyperplanes form two conjugacy classes of $11$-complements in $S$. In this case ${\rm P\Gamma L}_3(5)=S$ so the only almost simple group $G$ arising is the extension $G=A=S\rtimes{\rm C}_2$ of $S$ by its graph automorphism, transposing the actions of $S$ on points and hyperplanes; this has no transitive actions of degree a power of $11$, and thus no $11$-complements.

\medskip

Motivated by the Nagell--Ljunggren Conjecture (see Section~\ref{sec:simple}), we conjecture that this example is the only case to consider. Nevertheless, in case our conjecture is wrong, we should consider all other possibilities. If $G\le{\rm P\Gamma L}_d(q)$ and $G$ has $p$-complements, then so has $S$. The converse is true if and only if $|G:S|$ is coprime to $p$, so it is sufficient to consider whether or not there exist subgroups $G\le{\rm P\Gamma L}_d(q)$ with $|G:S|$ divisible by $p$. There are inclusions $S={\rm PSL}_d(q)\le {\rm PGL}_d(q)\le{\rm P\Gamma L}_d(q)$ with indices $(d,q-1)$ and $e:=|{\rm Gal}\,{\mathbb F}_q|$,
and Corollary~\ref{cor:d,q-1} shows that $p$ does not divide $(d,q-1)$, so we have:

\begin{prop}
Suppose that ${\rm PSL}_d(q)\le G\le {\rm P\Gamma L}_d(q)$ where these groups have natural degree $n=p^k$ for some prime $p$.
Then the point stabilisers in $G$ are $p$-complements if and only if the image $G/G\cap{\rm PGL}_d(q)$ of $G$ in ${\rm P\Gamma L}_d(q)/{\rm PGL}_d(q)\cong {\rm Gal}\,{\mathbb F}_q\cong{\rm C}_e$ has order coprime to $p$.
In particular, this holds for all $G$ such that ${\rm PSL}_d(q)\le G\le {\rm PGL}_d(q)$, and for all $G$ such that ${\rm PSL}_d(q)\le G\le {\rm P\Gamma L}_d(q)$ if $e$ is coprime to $p$. \hfill$\square$
\end{prop}


To summarise, in case~(B) the almost simple groups with socle ${\rm PSL}_d(q)$ which have $p$-complements are as follows:
\begin{itemize}
\item for $d=2$ the groups $S$ and $G$ in Proposition~\ref{prop:ASd=2};
\item for $d>2$ and $k=1$ any groups $G$ such that ${\rm PSL}_d(q)\le G\le {\rm P\Gamma L}_d(q)$ for which $n=(q^d-1)/(q-1)=p$;
\item for $d>2$ and $k>1$ any groups $G$ such that ${\rm PSL}_d(q)\le G\le {\rm P\Gamma L}_d(q)$ for which $n=(q^d-1)/(q-1)=p^k$, and
$|G:G\cap{\rm PGL}_d(q)|$ is coprime to $p$.
\end{itemize}
These groups have one conjugacy class of $p$-complements when $d=2$, and two classes when $d>2$. Only finitely many such groups are known in the first two cases, though it is conjectured that in each case there are infinitely many. The only known example of the third type is ${\rm PSL}_5(3)$,with $p=11$ and $k=2$, and it is conjectured that there are no others.

This, together with Proposition~\ref{prop:transnotB}, completes the proof of Theorem~\ref{th:mainpcomp}: note that the corresponding permutation groups are all primitive, so the $p$-complements listed there are all maximal subgroups.


\section{Extension to all primitive groups of prime power degree}\label{sec:allprim}

Having concentrated almost exclusively on almost simple groups, we now widen our scope to include {\sl all\/} primitive permutation groups of prime power degree. As noted by Cai and Zhang~\cite{CZ}, such a group $G$ of degree $p^k$ must have one of the following three O'Nan--Scott types (see~\cite[Chapter~4]{DM}):
\begin{itemize}
\item[(AS)] almost simple type, where the socle of $G$ is a nonabelian finite simple group $S$, and $S\le G\le{\rm Aut}\,S$;
\item[(HA)] holomorph affine type, where the socle of $G$ is an elementary abelian group $V$ of order $p^k$, regarded as a $k$-dimensional vector space over ${\mathbb F}_p$, and $G=V\rtimes G_0\le{\rm AGL}_k(p)$ where $G_0$ is a subgroup of ${\rm GL}_k(p)$ acting irreducibly on $V$;
\item[(PA)] product action type, where $G$ has socle $T^d$ for some nonabelian simple group $T$ and divisor $d>1$ of $k$, and $G\le U\wr K$ where $U$ is an almost simple primitive group of degree $p^e$ ($e=k/d$) with socle $T$, and $K$ is a transitive permutation group of degree $d=k/e$. (Here $U$ acts on a set $\Gamma$ of order $p^e$, and $G$ inherits the product action of $U\wr K$ on $\Omega=\Gamma^d$; the possibilities for $U$ are listed in Theorem~\ref{th:mainppower}, excluding the imprimitive groups in class~(A2).)
\end{itemize}
For the other types, the degree of $G$ cannot be a prime power. The structure of the socle, a property of the abstract group, determines the type AS, HA or PA of any faithful primitive representation of $G$, so they are all of the same type.

We can use this trichotomy to describe those groups with core-free maximal $p$-complements, since these are the primitive groups of degree $p^k$ in which the point stabilisers are $p'$-groups. 
\begin{itemize}
\item[(AS)] The almost simple groups with $p$-complements are listed in Theorem~\ref{th:mainpcomp}; the $p$-complements are all maximal since the corresponding permutation groups are primitive.
\item[(HA)] A group $G=V\rtimes G_0$ of holomorph affine type has $p$-complements if and only if $G_0$ is a $p'$-group, in which case by the Schur--Zassenhaus Theorem they are the conjugates of $G_0$;  they are maximal if and only if $G_0$ acts irreducibly on $V$.
\item[(PA)] When $G$ has product action type, the point stabilisers in $G$ are $p$-complements if and only if they are $p'$-groups, and for this it is necessary and sufficient that $K$ is a $p'$-group and $U$ has $p$-complements; such groups $U$ are listed in Theorem~\ref{th:mainpcomp}.
\end{itemize}


\section{Proof of Corollary~\ref{cor:regsgp}}\label{sec:regular}

By Proposition~\ref{prop:permgps}, in any transitive group in which the point stabilisers are $p$-complements (as for those in Theorem~\ref{th:mainpcomp}, for example) the Sylow $p$-subgroups act regularly. In fact, {\sl all\/} of the almost simple groups of prime power degree, as listed in Theorem~\ref{th:mainppower}, have regular subgroups, though they need not always be Sylow subgroups. The argument, for the various cases in that theorem, is as follows:
\begin{itemize}
\item[\rm(A1)] Cayley's Theorem states that every group $R$ of order $n$ is a regular subgroup of ${\rm S}_n$, so in particular this applies to each prime power $n=p^k\ge 5$. In this case we have $R\le{\rm A}_n$ unless $R\cong{\rm C}_n$ with $n=2^k$, so each group ${\rm A}_n$ in this class also has regular subgroups.
\item[\rm(A2)] Let $m=n/2=2^{k-1}$. Then ${\rm S}_m$ acts on $\Omega=\{\pm 1, \dots, \pm m\}$ as the product of its natural transitive actions on $\{1,\ldots, m\}$ and $\{\pm 1\}$. The $m$-cycle $a=(1,\ldots, m)\in{\rm S}_m$, which is odd, has cycles $(+1,-2,+3,\ldots, -m)$ and $(-1,+2,-3,\ldots,+m)$ on $\Omega$, while the involution $b=(1,m)(2,m-1)\ldots(m/2,m/2+1)$, which is even, inverts and transposes them, so $a$ and $b$ generate a dihedral group acting regularly on $\Omega$.
\item[\rm(B, C)] Apart from some possible exceptions in (B5), these groups all have $p$-complements, so their Sylow $p$-subgroups act regularly. In the exceptional cases, it is sufficient to apply this argument to the subgroup ${\rm PSL}_d(q)$.
\item[\rm(D)] These two groups have the same Sylow $3$-subgroups $P\cong{\rm C}_3\wr{\rm C}_3$ of order $3^4$. In addition to the base group ${\rm C}_3^3$, which is intransitive, $P$ has three regular subgroups $R$, which are non-abelian groups of order $3^3$ and exponent $3^2$.
\end{itemize}

Groups of holomorph affine type have regular normal subgroups, consisting of translations. In a primitive group $G$ of product action type, as described in Section~\ref{sec:allprim}, the simple group $T$ has a subgroup $R$ acting regularly on $\Gamma$, by the above argument, so $R^d$ is a subgroup of $G$ acting regularly on $\Omega=\Gamma^d$. This completes the proof of Corollary~\ref{cor:regsgp}, that every primitive permutation group of prime power degree has a regular subgroup. \hfill$\square$

\medskip

Corollary~\ref{cor:regsgp} does not, in general, extend to primitive groups for which the degree is not a prime power.

\medskip

\noindent{\bf Example 11.1} Let $G={\rm PSL}_2(q)$, acting naturally with degree $q+1$ for some prime power $q\equiv 1$ mod~$(4)$. This action is doubly transitive, and hence primitive. Each point stabiliser has even order $q(q-1)/2$, so it contains an involution. Since involutions in $G$ are mutually conjugate, they all have fixed points. Any regular subgroup, having even order $q+1$, would contain an involution with no fixed points, which is impossible. Thus $G$ has no regular subgroups.

\medskip

Similarly, Corollary~\ref{cor:regsgp} does not extend to imprimitive\ groups of prime power degree, though in this case counter-examples are not quite so easy to find.
First let us define a permutation group to be {\sl minimal transitive\/} if it is transitive, but all its proper subgroups (equivalently, all its maximal subgroups) are intransitive. Clearly, a minimal transitive group which is not regular cannot have a regular subgroup. Since the Sylow $p$-subgroups of a transitive group of $p$-power degree are also transitive, it follows that a minimal transitive group of $p$-power degree must be a $p$-group. Since the maximal subgroups $M$ of a $p$-group $P$ are normal, of index $p$, the point stabilisers $P_{\alpha}$ must satisfy $P_{\alpha}\le M$ or $MP_{\alpha}=P$ ; however, the latter implies that $M$ is transitive, against our assumption, so $P_{\alpha}$ is contained in the intersection of all the maximal subgroups $M$ of $P$, that is, in the Frattini subgroup $\Phi(P)$ of $P$. Conversely, if $P$ is transitive and $P_{\alpha}\le\Phi(P)$, then every maximal subgroup of $P$ is intransitive, so $P$ is minimal transitive. To find groups of prime power degree with no regular subgroups, it is therefore sufficient to find $p$-groups $P$ which act faithfully and non-regularly with point stabilisers contained in $\Phi(P)$.

Transitive groups of prime degree clearly have regular subgroups, and by the above argument the same applies to transitive groups of degree $p^2$: a minimal transitive group $P$ satisfies $P_{\alpha}\le\Phi(P)$, so that $\Phi(P)$ has index $p$ or $p^2$ in $P$; in the first case $P$ is cyclic and hence regular, and in the second case $P_{\alpha}=\Phi(P)$, which is normal in $P$, so $P$ is again regular. However, this does not extend to higher powers of $p$.

\medskip

\noindent{\bf Example 11.2} For any prime $p>2$ and integer $e\ge 1$ let
\[P=\langle a, b, c, d\mid a^p=b^p=c^p=d^{p^e}=[a,b]=[b,c]=[c,a]=1,\, a^d=ab,\, b^d=bc,\, c^d=c\rangle.\]
This is a semidirect product of an elementary abelian normal subgroup $E=\langle a, b, c\rangle$ of rank $3$ by $D=\langle d\rangle\cong{\rm C}_{p^e}$, with $\Phi(P)=\langle [a,d]=b,\, [b,d]=c,\, d^p\rangle\cong{\rm C}_p\times{\rm C}_p\times{\rm C}_{p^{e-1}}$ and $b\not\in Z(P)=\langle c, d^p\rangle$, so $P$ acts faithfully on the cosets of $\langle b\rangle$ as a minimal transitive group of degree $p^k$ where $k=e+2\ge 3$. (This generalises an example in~\cite{Jon72} with $k=3$.) When $p=2$ a similar construction works for all $e\ge 2$; the case $e=1$ must be avoided since conjugation by $d$ induces an automorphism of $E$ of order $p^2$ when $p=2$, rather than $p$ when $p>2$. It follows that $P$, and indeed any group of degree $p^k$ containing $P$ as a Sylow $p$-subgroup, cannot have a regular subgroup, and must therefore by imprimitive by Corollary~\ref{cor:regsgp}.

\medskip

\noindent{\bf Example 11.3} For any prime $p$ and integer $e\ge 2$ let
\[P=\langle a, b\mid a^{p^e}=b^{p^2}=1, a^b=a^{p^{e-1}+1}\rangle=\langle a\rangle\rtimes\langle b\rangle\cong{\rm C}_{p^e}\rtimes{\rm C}_{p^2},\]
a metacyclic group of order $p^{e+2}$. (Note that $a\mapsto a^{p^{e-1}+1}$ extends to an automorphism of $\langle a\rangle$ of order $p$.) This has Frattini subgroup $\Phi(P)=\langle a^p, b^p\rangle$ of order $p^e$, and centre $Z(P)=\langle a^p\rangle$ of order $p^{e-1}$. Since $\langle b\rangle$ has trivial intersection with $Z(P)$, it follows that $P$ acts faithfully, with degree $p^{e+f}$, on the cosets of $\langle b^{p^f}\rangle$ for each $f=1, \ldots e-2$. Since $b^{p^f}\in\Phi(P)$ for each such $f$ we obtain minimal transitive permutation representations of $P$ of degree $p^k$, on the cosets of $\langle b^{p^f}\rangle$, for each $k=e+f=e+1,\ldots, 2e-2$. These representations are all imprimitive, and by its minimality $P$ has no regular subgroups in any of them. The same conclusion will apply to any transitive group of degree $p^k$ which has $P$, acting in this way, as a Sylow $p$-subgroup; the holomorph ${\rm Hol}\,{\rm C}_{p^e}\cong P\rtimes{\rm C}_{p-1}$ is an example.


\section{Proof of Corollary~\ref{cor:autequiv}}\label{sec:ineq}

Among the almost simple transitive groups of $p$-power degree, as listed in Theorem~\ref{th:mainppower}, the only groups with two inequivalent representations of the same degree $p^k$ are those in classes (B4), together with ${\rm PSL}_2(11)$ in class (C), with $k=1$, and those in class (B5), with $k>1$. While there are many (possibly infinitely many) such examples with $k=1$, only one is known with $k>1$, namely ${\rm PSL}_5(3)$ of degree $11^2$. However, for each $k>1$ there are further examples which are of product action type, rather than almost simple type.

\medskip

\noindent{\bf Example 12.1} Let $U$ be an almost simple group with ${\rm PSL}_d(q)\le U\le{\rm P\Gamma L}_d(q)$ in class (B4), having two representations of projective prime degree $p=(q^d-1)/(q-1)$. If $K$ is any transitive permutation group of degree $k$ then the wreath product $G=U\wr K=U^k\rtimes K$ has two inequivalent primitive representations of degree $p^k$. These are its product actions on the $k$th cartesian power $\Omega=\Gamma^k$ of the set $\Gamma$ of points or hyperplanes of ${\mathbb P}^{d-1}(q)$; these actions are primitive by~\cite[Lemma~2.7A]{DM}, since $U$ itself acts primitively and $K$ acts transitively. Having non-conjugate point stabilisers $H_i\wr K$ ($i=1, 2$), where $H_1$ and $H_2$ are non-conjugate stabilisers in $U$, these representations of $G$ are inequivalent, but they are transposed by an obvious outer automorphism of $G$, acting diagonally on $U^k$ with $H_1$ and $H_2$ transposed in each factor, and fixing $K$. Similar constructions with $U={\rm PSL}_2(11)$ and ${\rm PSL}_5(3)$ give groups $G=U\wr K$, each with two primitive representations of degree $11^k$ and $11^{2k}$ respectively, transposed by automorphisms of $G$. If $K$ is chosen to be a $p'$-group (and therefore $k$ is coprime to $p$), then in all cases the point stabilisers in $G$ are $p'$-groups (since each $H_i$ is), so they are $p$-complements in $G$.

\medskip

There are also examples of holomorph affine type with $k>1$, in this case having an unbounded number of inequivalent representations.

\medskip

\noindent{\bf Example 12.2} Let $G={\rm ASL}_2(q)=V\rtimes G_0$, where $V$ is the natural module ${\mathbb F}_q^2$ for $G_0={\rm SL}_2(q)$. The conjugacy classes of complements for $V$ in $G$ are in one-to-one correspondence with the elements of the first cohomology group $H^1(G_0,V)$ (see~\cite[Satz I.17.3]{Hup}, for example). Now Bell~\cite[Table~1]{Bell} has shown that if $q=2^e>2$ then $H^1({\rm SL}_2(q),V)$ has dimension~$1$ over ${\mathbb F}_q$, so there are $q=2^e$ conjugacy classes of complements for $V$ in $G$. These are the point stabilisers in $2^e$ mutually inequivalent faithful representations of $G$, all primitive and of degree $q^2=2^{2e}$. Of course, these complements are all isomorphic to $G/V$, and hence to ${\rm SL}_2(q)$; they are not $2$-complements, since they have even order. Similarly, it follows from~\cite{Bell} that there are two conjugacy classes of complements for the translation subgroup in ${\rm ASL}_3(2)={\rm AGL}_3(2)$, and hence two primitive representations of degree $8$; however, other groups ${\rm ASL}_n(q)$ have just one conjugacy class of complements, since $H^1({\rm SL}_n(q),V)$ is trivial~\cite{Bell}. In all cases where there is a choice of complements, the induced permutation group is isomorphic to ${\rm ASL}_n(q)$, but preserving a different affine structure on the permuted set. In other words, the corresponding permutation representations are mutually equivalent under ${\rm Aut}\,G$.

\medskip

Corollary~\ref{cor:autequiv} asserts that this behaviour, of inequivalent representations being mutually equivalent under automorphisms, applies to all groups $G$ with faithful primitive representations of the same prime power degree. To prove this, Theorem~\ref{th:mainppower} deals, by inspection, with cases where $G$ (together, therefore, with its faithful primitive permutation representations) is of type AS. (However, note that the group ${\rm PSL}_3(2)\cong{\rm PSL}_2(7)$ has primitive representations of two {\sl distinct\/} prime power degrees, namely $7$ and $8$ in classes (B4) and (B1).)
In the product action case, as described in Section~\ref{sec:allprim}, the group $G$ uniquely determines its socle $T^d$, and hence the simple group $T$ and the exponent $d$, and then the degree $|\Omega|=p^k$ of $G$ uniquely determines the degree $|\Gamma|=p^e$ ($e=k/d$) of $T$. By Theorem~\ref{th:mainppower} any two representations of degree $p^e$ of $T$ differ by an automorphism of $T$, so as in Example~12.1 the same applies to the corresponding product actions of $G$.
In the holomorph affine case, the point stabilisers $H$ in any primitive representation of degree $p^k$ are isomorphic to $G/V$, where $V$ is the socle of $G$, so they are isomorphic to each other; they all act on $V$ in the same way since their corresponding elements differ by elements of $V$, which is abelian. Thus the various permutation groups $V\rtimes H$ arising from different conjugacy classes of complements $H$ for $V$ are mutually isomorphic as permutation groups, so these conjugacy classes, and hence the associated permutation representations of $G$, are equivalent under automorphism of $G$, as in Example 12.2.  This proves Corollary~\ref{cor:autequiv}. \hfill$\square$

\medskip

The following examples show that in Corollary~\ref{cor:autequiv} one cannot relax the conditions that the representations should be faithful, primitive and of the same prime power degree. The first example deals with faithfulness.

\medskip

\noindent{\bf Example 12.3} The primitive groups $G^*={\rm AGL}_d(q)$ of prime power degree in Example~3.1 have mutually non-isomorphic. These are not core-free since they contain the normal translation subgroup $V$, so the permutation representations are not faithful. Thus the faithfulness condition cannot be omitted from Corollary~\ref{cor:autequiv}. 

\medskip

The following examples show that the same applies to primitivity.

\medskip

\noindent{\bf Example 12.4} For any prime $p>2$ and integer $e\ge 1$ let
\[P=\langle a, b, c, d\mid a^p=b^p=c^p=d^{p^e}=[a,b]=[b,c]=[c,a]=1,\, a^d=ab,\, b^d=bc,\, c^d=c\rangle,\]
as in Example~11.2. In addition to the faithful representation of degree $p^{e+2}$ on the cosets of $\langle b\rangle$ described there, which is minimal transitive since $b\in\Phi(P)=\langle b, c, d^p\rangle$, $P$ has a faithful representation of the same degree on the cosets of $\langle a\rangle$, and this is not minimal transitive since $a\not\in\Phi(P)$. These two representations are imprimitive, since the point stabilisers are proper subgroups of $E$, and they are inequivalent under ${\rm Aut}\,P$ since $\Phi(P)$ is a characteristic subgroup of $P$.

\medskip

In fact, groups can have an unbounded number of faithful, transitive but imprimitive representations of the same prime power degree, which are mutually inequivalent under automorphisms.

  \medskip
 
 \noindent{\bf Example 12.5} For any prime $p$ and integer $e\ge 1$ let $P={\rm C}_p\wr {\rm C}_{p^e}$.
One can regard the base group $B$ as the $p^e$-th cartesian power of ${\mathbb F}_p$, with the complement ${\rm C}_{p^e}$ cyclically permuting coordinates. For each $j=0, 1,\ldots, e$ let $H_j=\langle h\rangle\cong{\rm C}_p$ for some $h=(h_i)\in B$ with coordinates $h_1, h_2,\ldots$ which repeat with period $p^j$. Then each $H_j$ has $p^j$ conjugates in $P$, and in particular $H_0$ is the centre of $P$. The representations of $P$ on the $p^{p^e+e-1}$ cosets of $H_j$ for $j=1,\ldots, e$ are faithful, transitive, imprimitive and mutually inequivalent under ${\rm Aut}\,P$.

\medskip

The condition of having prime power degree is also essential in Corollary~\ref{cor:autequiv}; an easy way of demonstrating this is to give examples of groups with non-isomorphic core-free maximal subgroups of the same index, which is not a prime power.

\medskip

\noindent{\bf Example 12.6} In ${\rm PSL}_2(13)$ there are both alternating and dihedral maximal subgroups of order $12$ and index $91$, and there are similar examples based on maximal subgroups ${\rm S}_4$ in ${\rm PSL}_2(23)$, and ${\rm A}_5$ in ${\rm PSL}_2(59)$ and ${\rm PSL}_2(61)$. 
For an infinite family of examples, the stabilisers of points and of isotropic lines in the symplectic group ${\rm PSp}_4(q)$ for odd $q$ are non-isomorphic
maximal subgroups of index $q^3+q^2+q+1$ and order $q^4(q-1)(q^2-1)/2$;
the smallest of these simple groups, ${\rm PSp}_4(3)$, is isomorphic as an abstract group to the primitive group $W(E_6)'$ of degree $27$ in class (D). See~\cite[Section~3.5.6]{Wil} for this family of examples; there are similar families involving the maximal parabolic subgroups in $G_2(q)$ for $q$ coprime to $3$~\cite[Table~4.1]{Wil} and in $F_4(q)$ for odd $q$~\cite[Theorem~4.4 (ii), (iii)]{Wil}.

\medskip

\noindent{\bf Example 12.7} The following further instances of non-isomorphic maximal subgroups of the same index in a simple group are taken from~\cite{OATLAS} and \cite{CCNPW}:
\begin{itemize}
\item ${\rm M}_{12}$, two isomorphism classes of order 192 and index 495, stabilisers of tetrads and of linked fours;
\item ${\rm PSL}_3(7)$, two isomorphism classes of order 72 and index 26\,068;
\item ${\rm M}_{23}$, two isomorphism classes of order 40\,320 and index 253, stabilisers of duads and hexads;
\item ${\rm HS}$, two isomorphism classes of order 40\,320 and index 1\,100;
\item ${\rm G}_2(4)$, two isomorphism classes of order 184\,320 and index 1\,365;
\item ${\rm McL}$, two isomorphism classes of order 58\,320 and index 15\,400;
\item ${\rm McL}$, three isomorphism classes of order 40\,320 and index 22\,275;
\item ${\rm O'N}$, two isomorphism classes of order 25\,920 and index 17\,778\,376;
\item ${\rm Th}$, two isomorphism classes of order 944\,784 and index 96\,049\,408\,000;
\item ${\rm Th}$, two isomorphism classes of order 12\,000 and index 7\,562\,161\,990\,656;
\item ${\rm Co_1}$, two isomorphism classes of order 60\,000 and index 69\,296\,280\,109\,056.
\end{itemize}
(This list does not claim to be comprehensive. The apparent example in~\cite{OATLAS} of maximal subgroups $3.{\rm Fi}_{24}$ and $2^2.^2{\rm E}_6(2):{\rm S}_3$ of ${\rm M}$ of order $1\, 836\, 779\, 512\, 410\, 596\, 494\, 540\, 800$ is incorrect: $3.{\rm Fi}_{24}$ does not have this order.)

\medskip

As a digression, Examples 12.6 and 12.7 raise some interesting questions about maximal subgroups of finite simple groups:

\medskip

\noindent{\bf Question} Is there an upper bound on the number of isomorphism classes of maximal subgroups of a given order in a finite simple group? The maximum number in the above list, attained by subgroups of order 40\,320 in McL, is three.

\medskip

\noindent{\bf Question} Are there instances of isomorphic maximal subgroups which are not conjugate in ${\rm Aut}\,G$? The brief search which gave rise to Example~12.6 yielded none.


\section{Conjectures from Number Theory}\label{sec:NT}

Here we discuss two number-theoretic conjectures which are relevant to our work.

\subsection{The Bateman--Horn Conjecture}

One of the main problems we have addressed in this paper is whether the polynomial $q^{d-1}+q^{d-2}+\cdots+q+1$ can take prime or prime power values, where $q$ is a prime power, and if so, whether it does so finitely or infinitely many times. In this context, Schinzel's Hypothesis H, the Bateman--Horn Conjecture, and the Nagell-Ljunggren Conjecture are all relevant.

If a polynomial $f\in{\mathbb Z}[t]$ is to take prime values $f(t)$ for infinitely many $t\in{\mathbb N}$, then is it clearly necessary that
\begin{enumerate}
\item $f$ has a positive leading coefficient;
\item $f$ is irreducible;
\item $f$ is not identically zero modulo any prime.
\end{enumerate}
In 1857 Bunyakovsky~\cite{Bun} conjectured that these conditions are also sufficient. This has been proved only in the case where $\deg(f)=1$: this is Dirichlet's Theorem on primes in an arithmetic progression $at+b$. No other case has been proved, not even such an apparently simple example as the Euler--Landau Problem of whether there are infinitely many primes of the form $t^2+1$ for $t\in{\mathbb N}$. In 1904 Dickson~\cite{Dic} conjectured that a finite set of polynomials $f_1,\ldots,f_k\in{\mathbb Z}[t]$, all of degree $1$, simultaneously take prime values $f_i(t)$ for infinitely many $t\in{\mathbb N}$ if and only if they each satisfy Bunyakovsky's conditions~(1) and~(2), and their product $f=f_1\ldots f_k$ satisfies~(3); examples include the Twin Primes Problem, with $f_i(t)=t$ and $ t+2$, and the Sophie Germain Primes Problem, with $f_i(t)=t$ and $2t+1$. In 1958 Schinzel's Hypothesis H asserted the same, but without the restriction that $\deg(f_i)=1$ for all $i$. This conjecture, proved only in the case of Dirichlet's Theorem, includes all of the above examples, together with ours: if $q$ is the $e$th power of a prime, then we would like to know whether, for any fixed $d$ and $e$, the polynomials
\begin{equation}\label{eq:BHC}
f_1(t)=t\quad\hbox{and}\quad f_2(t)=t^{(d-1)e}+t^{(d-2)e}+\cdots+t^e+1
\end{equation}
can be simultaneously prime for infinitely many $t\in{\mathbb N}$.

In 1923 Hardy and Littlewood~\cite{HL} examined a number of specific cases, and gave heuristic asymptotic estimates $E(x)$ for the number $Q(x)$ of $t\le x$ such that each $f_i(t)$ is prime. These estimates were generalised in 1962 by Bateman and Horn~\cite{BH} to apply to all instances of Schinzel's Hypothesis: their conjecture (the BHC), as recently improved by Li~\cite{Li}, is that if $f_1,\ldots, f_k$ satisfy the conditions of Schinzel's Hypothesis, then
\begin{equation}\label{eq:BHC}
Q(x)\sim E(x):= C\cdot\int_a^{\infty}\prod_{i=1}^k\frac{1}{\ln f_i(t)}\,dt\quad\hbox{as}\quad x\to\infty,
\end{equation}
where $a$ is chosen large enough for the integral to avoid singularities where some $f_i(t)=1$, and
\begin{equation}\label{eq:HL}
C=\prod_p\left(1-\frac{1}{p}\right)^{-k}\left(1-\frac{\omega_f(p)}{p}\right)
\end{equation}
where the product is over all primes $p$, and $\omega_f(p)$ is the number of distinct roots of $f$ mod~$(p)$. The heuristic explanation for this is that the integral in (\ref{eq:BHC}) is obtained by applying the Prime Number Theorem, in its more precise form
\[\pi(x)\sim {\rm Li}(x)=\int_2^x\frac{dt}{\ln t},\]
 to the polynomials $f_i(t)$, on the assumption that they take prime or composite values independently. This assumption is generally false, and the constant $C$ is a product of correction factors, one for each prime $p$, replacing the probability that $k$ randomly and independently  chosen elements of ${\mathbb Z}_p$ are all non-zero with the probability that a randomly chosen element of ${\mathbb Z}_p$ is not a root of $f$. The infinite product in (\ref{eq:HL}) converges (conditionally and slowly) to a limit $C>0$, while the integral diverges as $x\to\infty$, so if the BHC is true then $Q(x)\to\infty$ and hence there are infinitely many $t$ such that each $f_i(t)$ is prime.
 
 It is often straightforward to calculate $\omega_f(p)$, for instance using Quadratic Reciprocity if $\deg(f)=2$. A good approximation to the infinite product in (\ref{eq:HL}) can then be found by taking the partial product over all primes $p$ up to, say, $10^8$. Maple will compute the integral in (\ref{eq:BHC}) accurately and almost instantly, by numerical integration. (In the original paper, Bateman and Horn replaced each term $\ln f_i(t)$ with $\deg(f_i)\ln t$, in effect ignoring the coefficients of the polynomials $f_i$; this was, no doubt, in order to simplify the numerical integration, performed on slow and primitive computers with the programming done `from scratch', before the days of Maple; the resulting formula is asymptotically correct, but the estimates are less accurate than those obtained from Li's formula, given here.) In this way estimates $E(x)$ can be obtained for large $x$, and then compared with the actual numbers $Q(x)$, obtained by evaluating the polynomials $f_i$ at $t=1, 2, \ldots $ and then testing each $f_i(t)$ for primality by using the Rabin--Miller test within Maple; this test is probabilistic, but the probability of error is negligible, and in any case a few false or omitted primes would be insignificant in relation to the millions often found.

\begin{table}[htbp]
\begin{center}
\begin{tabular}{c|c|c|c}
$x$ & $Q(x)$ & $E(x)$ & $E(x)/Q(x)$ \\
\hline
$1 \cdot 10^{10}$ &  15\,801\,827 & $1.579642126 \times 10^7$  & 0.9996579044 \\
$2 \cdot 10^{10}$ &  29\,684\,763 & $2.968054227 \times 10^7$  & 0.9998578150 \\
$3 \cdot 10^{10}$ &  42\,963\,858 & $4.296235691 \times 10^7$  & 0.9999650617 \\
$4 \cdot 10^{10}$ &  55\,877\,571 & $5.587447496 \times 10^7$  & 0.9999445924 \\
$5 \cdot 10^{10}$ &  68\,522\,804 & $6.852175590 \times 10^7$  & 0.9999847043 \\
$6 \cdot 10^{10}$ &  80\,962\,422 & $8.096382889 \times 10^7$  & 1.0000173771 \\
$7 \cdot 10^{10}$ &  93\,236\,613 & $9.323905289 \times 10^7$  & 1.0000261688 \\
$8 \cdot 10^{10}$ & 105\,372\,725 & $1.053741048 \times 10^8$  & 1.0000130940 \\
$9 \cdot 10^{10}$ & 117\,383\,505 & $1.173885689 \times 10^8$  & 1.0000431394 \\
$10^{11}$         & 129\,294\,308 & $1.292974079 \times 10^8$  & 1.0000239757 
\end{tabular}
\end{center}
\vspace{2mm}
\caption{The second column gives the numbers $Q(x)$ of primes 
$p=q^2+q+1$ for primes $q\le x=i\cdot 10^{10}$, where $i=1,\ldots,10$.
The third column gives the corresponding Bateman--Horn estimates $E(x)$ for $Q(x)$,
and the fourth column gives the ratios $E(x)/Q(x)$.}
\label{tab:n=3BHratios}
\end{table}

The authors of~\cite{JZ2}, looking for groups ${\rm PSL}_d(q)$ of prime degree, avoided the difficult and well-trodden territory of Fermat and Mersenne primes arising when $d=2$, and instead chose the simplest remaining case, where $d=3$ and $e=1$. Thus $f_1(t)=t$ and $f_2(t)=t^2+t+1$, and we are counting primes $q\le x$ such that ${\rm PSL}_3(q)$ has prime degree $q^2+q+1$. The resulting estimates $E(x)$ and actual numbers $Q(x)$ are shown in Table~\ref{tab:n=3BHratios}.
 
Other small values of $d$ and $e$ were also tested in~\cite{JZ2}. The results were similar, but not so persuasive, as the values of $f_2(t)$ increase so rapidly that there are relatively few primes within reasonable computing range. This led the authors of~\cite{JZ2} to conjecture that for each prime $d\ge 3$ and each integer $e\ge 1$ there are infinitely many groups ${\rm PSL}_d(q)$ of prime degree (and thus possessing $p$-complements) where $q$ is the $e$th power of a prime. Equally convincing evidence for the BHC has been found in problems arising from other areas, such as block designs~\cite{JZ3} and orders of simple groups~\cite{JZ4} for example.

\subsection{The Nagell--Ljunggren Conjecture}

The Nagell--Ljunggren equation
\[
\frac{x^d-1}{x-1}=y^k
\] 
originated in the papers~\cite{Nag1, Nag2} by Nagell and~\cite{Lju} by Ljunggren. In~\cite{BM} Bugeaud and Mignotte have conjectured that the only solutions $x, y\in{\mathbb Z}$ with $|x|$, $|y|$, $k>1$ and $d>2$ are the four given in (\ref{eq:NLknown}). For instance, Ljunggren proved that there are no further solutions with $k=2$.
There is an excellent survey of the situation up to 2002 in~\cite{BM}, and Bennett and Levin have given a very useful and more recent update in~\cite{BL}. Combining results from a number of papers, they show that (changing their standard notation to that used here) if $(x,y,d,k)$ is a solution, other than one of those currently known, then
\begin{enumerate}
\item $k$ is odd;
\item $d$ is not divisible by $2$ or $3$;
\item $|x|\ge 10^4$, and $x$ has a prime divisor $r\equiv 1$ mod~$(k)$.
\end{enumerate}
The main result of~\cite{BL} is that $d$ is a product of at most three primes, counted with multiplicity. Progress on this general problem seems to be very difficult to achieve, and a complete solution is not in sight. A more modest conjecture in~\cite{BM}, that the equation has only finitely many solutions, seems more attainable: for instance, Bugeaud and Mignotte show that this would follow from a proof of the $abc$ conjecture.

However, in our group-theoretic application we are in a more special situation, considering the analogous equation
\begin{equation}\label{eq:NLqp}
\frac{q^d-1}{q-1}=p^k
\end{equation}
with $q$ now a prime power and $p$ prime. We conjecture that with these conditions, the only solution with $d>2$ and $k>1$ is $(3^5-1)/(3-1)=11^2$. For instance, in this case we have a fairly simple proof of a much stronger result, namely that $d$ is prime (see Corollary~\ref{cor:dprime}(a)), and hence $d\ge 5$ by (2). Similarly, (3) implies that any new solution has $q=r^e$ for a prime $r\equiv1$ mod~$(k)$, so $r\equiv 1$ mod~$(2k)$ by (1) and hence $r\ge 7$. A weaker conjecture, that equation~(\ref{eq:NLqp}) has only finitely many solutions, may be more achievable: for instance, it follows from~\cite[Theorem~5]{BM} that there are only finitely many if either $d$, $q$ or $p$ is fixed. Perhaps the powerful techniques developed to address the general Nagell-Ljungrren Conjecture may, in the hands of expert number-theorists, allow further progress towards solving this more restricted problem.

\medskip

Finally, we emphasise the remarkable contrast between the apparent abundance of solutions for equation~(\ref{eq:NLqp}) when $k=1$ (see Table~\ref{tab:n=3BHratios} for $d=3$), and their rarity when $k>1$ (see above for $d\ge 3$, and Proposition~\ref{prop:d=2} for $d=2$). An explanation of this phenomenon would be very welcome.


\bigskip

\noindent{\bf Acknowledgments} The authors are grateful to Danila Revin and Rob Wilson for some helpful comments on $p$-complements and on maximal subgroups of finite simple groups, and to Alexander Zvonkin for the computations presented in Table~\ref{tab:n=3BHratios}.


\bigskip

School of Mathematical Sciences, University of Southampton, Southampton SO17 1BJ, UK

G.A.Jones@maths.soton.ac.uk

\medskip

Department of Mathematics, Istanbul Okan University, Istanbul, Turkey

sezgin.sezer@okan.edu.tr
	
\end{document}